\documentclass{amsart}

\usepackage{latexsym}
\usepackage{amsfonts}
\newtheorem{theorem}{Theorem}[section]

\newtheorem{proposition}[theorem]{Proposition}
\newtheorem{remark}[theorem]{Remark}

\theoremstyle{definition}

\theoremstyle{remark}

\theoremstyle{remark}

\begin{document}

\newcommand\Pa{Painlev\'e}
\newcommand\Pt{{\rm P}_{\rm{\scriptstyle II}}}
\newcommand\Ptn{{\rm P}^{(n)}_{\rm{\scriptstyle II}}}
\newcommand\Lr{{\mathcal L}}
\newcommand\Pol{{\mathcal P}}
\newcommand\Se{{\mathcal S}}
\newcommand\A{{\mathcal A}}
\newcommand\M{{\mathsf M}}
\newcommand\Complex{{\mathbb C}}
\newcommand\Real{{\mathbb R}}
\newcommand{\ID}{\mathbb{1}}
\newcommand\Integer{{\mathbb Z}}
\def\ddt{{{\rm d}\over{\rm d} t}}
\def\ddx{{{\rm d}\over{\rm d} x}}
\def\ddtt{{{\rm d}\over{\rm d}\tau}}
\def\ddz{{{\rm d}\over{\rm d} z}}
\def\thi{{\vartheta_\infty}}
\def\eps{\varepsilon}
\def\arg{{\rm arg}}
\def\iff{\Longleftrightarrow}
\def\U{{\mathcal U}}

\title[Tritronqu\'ee Solutions of the Second Painlev\'e Hierarchy]
{Existence and
Uniqueness of Tri-tronqu\'ee Solutions of the second Painlev\'e
hierarchy}

\author{N. Joshi}
\address{School of Mathematics and Statistics F07, University of Sydney,
NSW2006 Sydney, Australia}
\email{nalini@maths.usyd.edu.au}
\thanks{Research supported by Australian Research Council Fellowship
\#F69700172, Discovery Grant \#DP0208430, and by Engineering and
Physical Sciences Research Council Fellowship \#GR/M28903.}
\author{M. Mazzocco}
\address{Mathematical Institute, University of Oxford, 24-29 St Giles, Oxford
OX1 3LB, UK and DPMMS, Wilberforce Road, Cambridge
CB3 0WB, UK.}
\email{M.Mazzocco@dpmms.cam.ac.uk}

\begin{abstract}
The first five classical Painlev\'e equations are known to have solutions
described by divergent asymptotic power series near infinity. Here we
prove that such solutions also exist for the infinite
hierarchy of equations associated with the second Painlev\'e equation.
Moreover we prove that these are unique in certain sectors near infinity.
\end{abstract}
\maketitle

\section{Introduction}
The {\it second Painlev\'e hierarchy}\/ is an infinite sequence of
nonlinear ordinary differential equations containing
\begin{equation}
\label{p2}
\Pt :\quad\quad V'' = 2V^3+x\,V+\alpha ,
\qquad V=V(x),\ \alpha\ \rm{const},
\end{equation}
as its simplest equation. This hierarchy is a
symmetry reduction of the mKdV hierarchy \cite{AS,Ai,FN},
and for this reason, it is believed that its elements
all posses the {\em Painlev\'e Property}\/, i.e.
all the movable singularities of all solutions are poles
\cite{Pain}. We conjecture that, in fact, all solutions of every equation
in the hierarchy are meromorphic in the complex
$x$-plane.

Another remarkable property of the equations belonging to
the second Painlev\'e hierarchy, denoted by $\Ptn$, is
their {\it irreducibility,}. That is, for generic values of the parameter
$\alpha_n$ in each equation of the hierarchy,
there is no transformation, within a certain class
described in \cite{Um}, that maps any of these
equations to a linear equation or to a lower order nonlinear ordinary
differential equation. This result has been proved only for the case of
$\Pt$ (see \cite{UW}), but
we conjecture that it is the case for all the other equations
in the hierarchy. In fact, we show evidence that no $n$-th
order member of the hierarchy can be reduced to a lower-order
member of the hierarchy.

The main aim of this paper is the asymptotic study of solutions of
each $n$-th equation $\Ptn$ in the second \Pa\ hierarchy.
Since $\infty$ is a non-Fuchsian singularity for each equation $\Ptn$,
one expects that the generic solution of $\Ptn$ possesses asymptotic
behaviours that are given by (meromorphic combinations of)
hyperelliptic functions (see \cite{Li}). We concentrate
here, however, on solutions described by divergent power series
as $x\to\infty$. In particular, we focus on
solutions that possess no poles outside a circle of sufficiently large
radius in certain sectors of the $x$-plane.  We
show that there exist unique true solutions with
such behaviour in certain sectors of the $x$-plane. These have not
been proved before for $n\ge 2$.

To be more precise,
consider $\Pt$.
The asymptotic study of its solutions, in the limit $x\to\infty$, began
with the work of Boutroux \cite{boutroux:I}.
The sectors of validity of its asymptotic behaviours
are described by six rays: $\arg(x)=j\pi/3$, $j=0, \ldots, 5$. All
solutions are meromorphic in the complex $x$-plane and the general
solution is asymptotic to Jacobian elliptic functions
whose modulus varies with $\arg(x)$ (see \cite{njmdk:conn1}).
Locally, the poles of the solution of $\Pt$
are aligned with the lattice of periodicity of its asymptotic behaviour.
However, as $\arg(x)$ changes within each sector of validity,
this lattice slowly changes.
Since the elliptic function has poles in each period parallelogram,
so does the solution of $\Pt$. In other words, the general
solution has an infinite
number of poles within each sector.

However, there also exist two types of one-parameter families of
solutions, called {\em
tronqu\'ee} by Boutroux, that possess no poles whatsoever in an
annular sector
\[
\Omega_{j}=\Bigl\{x\in\Complex\bigm | |x|>|x_{0}|,\
\bigl|\arg(x)-j\pi/3\bigr|<\pi/3\Bigr\},\qquad j=0,\dots,5,
\]
for some $x_{0}$. For a given $j$, as $|x|\to\infty$,
$x\in\Omega_{j}$, such a tronqu\'ee solution either has
asymptotic behaviour
\begin{equation}
 \label{yinfty}
y(x)=\left(\!\frac{-x}{2}\right)^{1/2}\left(1+{\mathcal O}
(x^{-\frac 32(1-\epsilon)})\right),
\end{equation}
for an appropriate choice of branch of the square root,
or
\begin{equation}
 \label{yzero}
y(x)=\left(\!\frac{-\alpha}{x}\right)\left(1+{\mathcal O}
(x^{-\frac 32(1-\epsilon)})\right),
\end{equation}
for some $\epsilon>0$. Boutroux characterised each sector $\Omega_{j}$ by
its ray of symmetry $\arg(x)=j\pi/3$. He showed that there also exist
{\em tri-tronqu\'ee} solutions
that are asymptotically pole-free along three successive such rays.
There are six such solutions. In this paper, we prove the existence
and uniqueness of the analogous {\em tri-tronqu\'ee} solutions
for all differential equations of the second Painlev\'e hierarchy.

The hierarchy we consider arises as a symmetry reduction
of the KdV (Korteweg-de Vries) hierarchy which is defined by
\begin{equation}
\label{kdv} \partial_{t_{n}} U+\partial_{\xi}\Lr_{n}\{U\} =0,
\ n\ge 1
\end{equation}
where
\begin{subequations}
\begin{equation}
\label{recursion} \partial_{\xi}\Lr_{n+1}\{U\}=
\bigl(\partial_{\xi\xi\xi}
         +4U\partial_{\xi}+2U_{\xi}\bigr)\Lr_{n}\{U\}
    \end{equation}
    \begin{equation}
\label{first} \Lr_{1}\{U\}=U
\end{equation}
    \end{subequations}
The fact that at each step equation (\ref{recursion}) can be
integrated to obtain the differential operator $\Lr_{n+1}$ was proved
in \cite{lax}.

After the transformation $U=W_{\xi}-W^{2}$, and reduction
$W(\xi)=V(x)/\bigl((2n+1)t_{2n+1}\bigr)^{1/(2n+1)}$,
$x=\xi/\bigl((2n+1)t_{2n+1}\bigr)^{1/(2n+1)}$ (see \cite{cjp} for
details), we get the $\Pt$ hierarchy
\begin{equation}
\label{p2hier}
\Ptn :\quad\quad \left(\frac{d}{dx}+2 V\right)\Lr_{n}\bigl\{V_{x}-V^{2}\bigr\}
=x V+\alpha_{n},\ n\ge 1
    \end{equation}
where $\alpha_{n}$ are constants and $\Lr_{n}$ is the operator defined by
Equation (\ref{recursion})
with $\xi$ replaced by $x$. We give a more explicit form of
(\ref{p2hier}) in Proposition \ref{nalini1} below. For $n=1$, Equation
(\ref{p2hier}) is $\Pt$, whereas for $n=2$ (noting that
$\Lr_{2}\{U\}=U_{\xi\xi}+3U^{2}$) it is
\begin{equation}
 \label{p24}
 V_{4x}-10 V^{2}V_{2x}-10V{V_{x}}^{2}+6V^{5}= x V+\alpha_{2},
 \end{equation}
where we use the notation
\[V_{mx}:=\frac{d^{m} V}{dx^{m}}.\]
Below we will also use $V_m\equiv V_{mx}$.

In this paper, we show that each equation $\Ptn$ of the
hierarchy possesses solutions
that are free of poles in sectors of angular width $2n\pi/(2n+1)$
as $|x|\to\infty$ and describe their asymptotic behaviour.
Furthermore, there exist unique solutions that are pole-free in
sectors of angular width $4 n\pi/(2n+1)$. We call such solutions
{\em tri-tronqu\'ee} solutions. 
In the second order case these solutions occur as important
solutions of transition phenomena for PDEs such as the modified
Korteweg de Vries equation (see \cite{Hab1}) and for ODEs with slowly
changing parameters (see \cite{Hab2}). We expect that the tri-tronqu\'ee
solutions of $\Ptn$
will occur in higher order PDEs with transition phenomena.

Our main result is

\begin{theorem}\label{main}
    For each integer $n\ge 1$, there exists $x_{0}$, $|x_{0}|>1$,
    and sectors
    \begin{eqnarray*}
      \Se_n&=&\left\{x\in\Complex\bigm| |x|>|x_{0}|,
        |\arg(x-x_{0})|<n\pi/(2n+1)\right\}\\
  \Sigma_n&=&\left\{x\in\Complex\bigm| |x|>|x_{0}|,
    |\arg(x-x_{0})|<2n\pi/(2n+1)\right\}
  \end{eqnarray*}
    in which the following results hold.
\begin{list}{\rm{\arabic{enumi}}.}{\usecounter{enumi}
\setlength{\leftmargin}{8truept}}
\item Equation (\ref{p2hier}) has formal solutions
    \begin{equation}
       \label{n2.5}
        V_{f, \infty}=\left(\frac{(-1)^{n}x}{2c_{n}}\right)^{\frac{1}{2n}}
        \sum_{k=0}^{\infty} a_{k}^{(n)}\left(\frac{2n}{2n+1}
        x^{\frac{2n+1}{2n}}\right)^{-k},\ a_{0}^{(n)}=1
        \end{equation}
    \begin{equation}
      \label{n3.5}
        V_{f, 0}=\left(\frac{-2\,\alpha_{n}}{(2n+1)\,x}\right)
        \sum_{k=0}^{\infty} b_{k}^{(n)}\left(\frac{2n}{2n+1}
        x^{\frac{2n+1}{2n}}\right)^{-k},\ b_{0}^{(n)}=1
        \end{equation}
    where
    \[c_{n}:=\frac{2^{2n-1}\Gamma(n+1/2)}{\Gamma(n+1)\Gamma(1/2)}\]
and the coefficients $a_{k}^{(n)}$, $b_{k}^{(n)}$, $k\ge 1$,
are given by substitution.
\item
   In each sector $\Se_n$ there exist true solutions
$V_{\infty}$ and $V_{0}$ of Equation (\ref{p2hier}) with asymptotic behaviour
    \begin{align*}
        V_{\infty}&\underset{{\overset{x\to\infty}{x\in\Se_n}}}{\sim}
          V_{f,\infty}\\
        V_{0}&\underset{{\overset{x\to\infty}{x\in\Se_n}}}{\sim} V_{f,0}
        \end{align*}
\item The true solutions $V_{\infty}$ and $V_{0}$ of
    Equation (\ref{p2hier}) are unique in proper sub-sectors of $\Sigma_n$
containing respectively $2n$ among the following half-lines
\begin{equation}
{\rm arg}(x) = {2j+1\over 2n+1}\pi,
\quad j=0,\dots,2n,
\label{n2.66}\end{equation}
or
\begin{equation}
{\rm arg}(x) = \pi + {2j+1\over 2n+1}\pi,
\quad j=0,\dots,2n.
\label{n3.66}\end{equation}
\end{list}
    \end{theorem}


The theorem above is proved in Section 2.

\begin{remark} We mention that although the elements of the $\Pt$
hierarchy are integrable (see \cite{CJM}), we make no use of this
fact in this paper.
To our knowledge there is no general result linking integrability to
existence of tri-tronqu\'ee-type solutions (for example in the case of
Painlev\'e VI equation such solutions have not yet been found).
\end{remark}


\section{Proof of Theorem \ref{main}}

We start the proof by deriving a more explicit expression for
the $\Pt$ hierarchy.
\begin{proposition}\label{nalini1}
The differential equations $\Ptn$ have the form
\begin{equation}
V_{2 n}= P_{2n-1}(V_0,V_1,\dots,V_{2n-2}) + x V_0 + \alpha_{n} +
\beta_{n} V_0^{2n+1},
\label{n1}\end{equation}
where $V_m:=\frac{d^{m} V}{dx^{m}}$, $P_{2n-1}$ is a polynomial
in $V_0,V_1,\dots,V_{2n-2}$ of degree $2n-1$, of the form
\begin{equation}\label{polyn}
P_{2n-1}=\underset{{\overset{
\scriptstyle \langle{\bf k}\rangle=2n+1}{\scriptstyle k_0\leq 2n-1}}}
{{\boldsymbol{\Sigma}}}
b_{k_0,\dots,k_{2n-2}} V_0^{k_0} V_1^{k_1}\dots
V_{2n-2}^{k_{2n-2}}.
\end{equation}
Here $k$ is a multi-index $k=(k_0,\ldots, k_{2n-2})$, with norm
\[\langle{\bf k}\rangle:=\sum_{p=0}^{2n-2}(p+1)k_p,\]
$b_{k_0,\dots,k_{2n-2}}$ are constants (some of which may be zero),
and
\[
\beta_{n}=(-1)^{n+1} 2^{2n}
{\Gamma(n+{1\over2})\over\Gamma(n+1)\Gamma({1\over2})}.
\]
\end{proposition}
\begin{proof}
The result (\ref{n1}) can be derived from (\ref{p2hier}) by proving that
\begin{equation}\begin{split}
\Lr_{n}(V_1-V_0^2) &= V_{2 n-1}+ \tilde\beta_n V_0^{2n}\\
&\quad +
\underset{{\overset{
\scriptstyle \langle{\bf k}\rangle=2n}{\scriptstyle k_0\leq 2n-2}}}
{{\boldsymbol{\Sigma}}}
a_{k_0,\dots,k_{2n-2}} V_0^{k_0}V_1^{k_1}\dots
V_{2n-2}^{k_{2n-2}}
\label{Ln2}
\end{split}
\end{equation}
where  $\tilde\beta_{n}=-{\beta_n\over2}$, and $a_{k_0,\dots,k_{2n-2}}$
denote constants.
Let us derive (\ref{n1}) from (\ref{Ln2}) first and then prove
(\ref{Ln2}) by induction.

The $n$-th equation $\Ptn$ (see Equation (\ref{p2hier})) has the form
\begin{equation*}\begin{split}
\Bigl({\partial\over\partial x}+2 V_0\Bigr)\Bigg(V_{2 n-1}&+
\tilde\beta_n V_0^{2n}+\underset{{\overset{
\scriptstyle \langle{\bf k}\rangle=2n}{\scriptstyle k_0\leq 2n-2}}}
{{\boldsymbol{\Sigma}}}
a_{k_0,\dots,k_{2n-2}} V_0^{k_0} V_1^{k_1}\dots
V_{2n-2}^{k_{2n-2}}\Bigg)=\\
&=x V_0+\alpha_n
\end{split}
\end{equation*}
that is
\begin{equation*}
\begin{split}
{V}_{2 n}&+ 2n\tilde\beta_n V_0^{2n-1}V_1+ \\
&+\underset{{\overset{
\scriptstyle \langle{\bf k}\rangle=2n}{\scriptstyle k_0\leq 2n-2}}}
{{\boldsymbol{\Sigma}}}
a_{k_0,\dots,k_{2n-2}}\Biggl(k_0 V_0^{k_0-1} V_1^{k_1+1}\dots
V_{2n-2}^{k_{2n-2}}+\\
&\qquad+k_1 V_0^{k_0} V_1^{k_1-1}V_2^{k_2+1}\dots
V_{2n-2}^{k_{2n-2}}+
\dots+\\
&\qquad+k_{2n-2} V_0^{k_0} V_1^{k_1}\dots
V_{2n-2}^{k_{2n-2}-1}V_{2n-1} \Biggr)+\\
&+ 2 V_0 V_{2 n-1}+2 \tilde\beta_n V_0^{2n+1}+\\
&+2\underset{{\overset{
\scriptstyle \langle{\bf k}\rangle=2n}{\scriptstyle k_0\leq 2n-2}}}
{{\boldsymbol{\Sigma}}}
a_{k_0,\dots,k_{2n-2}} V_0^{k_0+1} V_1^{k_1}\dots
V_{2n-2}^{k_{2n-2}}=\\
&\quad = x V_0+\alpha_n
\end{split}
\end{equation*}
$\Rightarrow$
\begin{equation*}
\begin{split}
V_{2 n}+ 2 \tilde\beta_n V_0^{2n+1}&-
\underset{{\overset{
\scriptstyle \langle{\bf k}\rangle=2n+1}{\scriptstyle k_0\leq 2n-1}}}
{{\boldsymbol{\Sigma}}}
b_{k_0,\dots,k_{2n-1}} V_0^{k_0} V_1^{k_1}\dots
V_{2n-1}^{k_{2n-1}}=\\
&\qquad =x V_0+\alpha_n
\end{split}
\end{equation*}
for some suitable constants $b_{k_0,\dots,k_{2n-1}}$. For
$\beta_n=-2\tilde\beta_n$ we obtain (\ref{n1}) as desired.

We now want to prove (\ref{Ln2}) by induction. It is trivially true
for $n=1$ with $\beta_1=2$. Let us assume it true for $n$ and prove
it for $n+1$. Observe that
$$
\Lr_{n+1} = {\partial^2\Lr_n\over\partial x^2} + 2(V_1-V_0^2)\Lr_n
+ 2 \int (V_1-V_0^2){\partial\Lr_n\over\partial x}{\rm d}x,
$$
where $(V_1-V_0^2){\partial\Lr_n\over\partial x}{\rm d}x$
is always an exact form (see \cite{lax}). We obtain
$$
\begin{array}{ll}
\Lr_{n+1} =& V_{2n+1}+2n(2n-1)\tilde\beta_n V_0^{2n-2}V_1^2+
2n\tilde\beta_n V_0^{2n-1}V_2 + \\
&+\underset{{\overset{
\scriptstyle \langle{\bf k}\rangle=2n}{\scriptstyle k_0\leq 2n-2}}}
{{\boldsymbol{\Sigma}}}
a_{k_0,\dots,k_{2n-2}}\big(k_0(k_0-1) V_0^{k_0-2} V_1^{k_1+2}\dots
V_{2n-2}^{k_{2n-2}}+\\
&+k_0 (k_1+1) V_0^{k_0-1} V_1^{k_1} V_2^{k_2+1}\dots
V_{2n-2}^{k_{2n-2}}+\dots+\\
&+ k_0 k_{2n-2} V_0^{k_0-1} V_1^{k_1+1} \dots
V_{2n-2}^{k_{2n-2}-1}V_{2n-1}+\dots+\\
&+k_{2n-2}k_0 V_0^{k_0-1} V_1^{k_1+1}\dots
V_{2n-2}^{k_{2n-2}-1}V_{2n-1}+\dots+\\
&+k_{2n-2}(k_{2n-2}-1) V_0^{k_0} V_1^{k_1}\dots
V_{2n-2}^{k_{2n-2}-2}V_{2n-1}^2+\\
&+k_{2n-2}V_0^{k_0} V_1^{k_1}\dots V_{2n-2}^{k_{2n-2}-1}V_{2n}
 \big)+ 2 V_1 V_{2n-1} + 2 V_1 \tilde\beta_n V_0^{2n} +\\
&+2 V_1
\underset{{\overset{
\scriptstyle \langle{\bf k}\rangle=2n}{\scriptstyle k_0\leq 2n-2}}}
{{\boldsymbol{\Sigma}}}
a_{k_0,\dots,k_{2n-2}} V_0^{k_0} V_1^{k_1}\dots
V_{2n-2}^{k_{2n-2}} - 2 V_0^2  V_{2n-1} -\\
&- 2\tilde\beta_n V_0^{2n+2}- 2 V_0^2
\underset{{\overset{
\scriptstyle \langle{\bf k}\rangle=2n}{\scriptstyle k_0\leq 2n-2}}}
{{\boldsymbol{\Sigma}}}
a_{k_0,\dots,k_{2n-2}} V_0^{k_0} V_1^{k_1}\dots
V_{2n-2}^{k_{2n-2}}+\\
&+2 \int (V_1-V_0^2)\bigg[V_{2n}+ 2n\tilde\beta_n V_0^{2n-1}V_1 +\\
&+\underset{{\overset{
\scriptstyle \langle{\bf k}\rangle=2n}{\scriptstyle k_0\leq 2n-2}}}
{{\boldsymbol{\Sigma}}}
a_{k_0,\dots,k_{2n-2}}\big(k_0 V_0^{k_0-1} V_1^{k_1+1}\dots
V_{2n-2}^{k_{2n-2}}+\\
&+k_1 V_0^{k_0} V_1^{k_1-1}V_2^{k_2+1}\dots V_{2n-2}^{k_{2n-2}}+
\dots+\\
&+k_{2n-2} V_0^{k_0} V_1^{k_1}\dots
V_{2n-2}^{k_{2n-2}-1}V_{2n-1} \big)
\bigg]{\rm d}x.\\
\end{array}
$$
Thus
$$
\begin{array}{ll}
\Lr_{n+1} =& V_{2(n+1)-1}- 2\tilde\beta_n V_0^{2(n+1)} -
{4n\over2n+2}\tilde\beta_n V_0^{2(n+1)} +\\
&+ \underset{{\overset{
\scriptstyle \langle{\bf k}\rangle=2n}{\scriptstyle k_0\leq 2(n+1)-2}}}
{{\boldsymbol{\Sigma}}}
\tilde a_{k_0,\dots,k_{2(n+1)-2}} V_0^{k_0} V_1^{k_1}\dots
V_{2(n+1)-2}^{k_{2(n+1)-2}}\\
\end{array}
$$
that has the required form for
$\tilde\beta_{n+1}=-2\tilde\beta_n{2n+1\over n+1}$, i.e.
$\beta_{n+1}=-2\beta_n{2n+1\over n+1}$. From $\beta_1=2$ we obtain
the right value of $\beta_n$.
\end{proof}

\begin{proposition}
    For each $n\ge 1$, the change of variables
 \begin{equation}
    \label{boutroux}
    V(x)=f(x)\, u(z), \qquad  z=g(x)
    \end{equation}
where
    \[ f(x)=x^{1/(2n)},\qquad g(x)=\frac{2n}{2n+1}\,x^{(2n+1)/2n} \]
maps $\Ptn$ to
\begin{equation}
\begin{array}{ll}
{{\rm d}^{2n}\over{\rm d} z^{2n}} u
=&-\sum_{l=0}^{2n-1} a_{2n,l}z^{l-2n}
{{\rm d}^l\over{\rm d} z^l} u(z)+{2n\over 2n+1} u(z) +
\beta_n u^{2n+1}+{2n\over 2n+1}{\alpha\over z}+\\
&+z^{-2n-1}\sum_{m_0,\dots,m_{2n-2}=0}^{2n-1}
a_{m_0,\dots,m_{2n-2}}(z)\prod_{l=0}^{2n-2}
\left({{\rm d}^l\over{\rm d} z^l} u(z)\right)^{m_l},\\
\end{array}
\label{bou4}\end{equation}
where $a_{m_0,\dots,m_{2n-2}}(z)$ are
polynomials of degree $\le 2n$ and $a_{2n,l}$ are some constant
coefficients and the sum $\sum_{m_0,\dots,m_{2n-2}=0}^{2n-2}$ is
zero for $n=1$. Furthermore, equation (\ref{bou4})
admits formal series expansions
\begin{equation}
\label{series1}
        u_{f, \infty}=\left(\frac{(-1)^{n}}{2c_{n}}\right)^{\frac{1}{2n}}
        \sum_{k=0}^{\infty} a_{k}^{(n)}z^{-k},\ a_{0}^{(n)}=1
        \end{equation}
\begin{equation}
     \label{series2}
u_{f, 0}=\left(\frac{-2\,\alpha_{n}}{2n\,z}\right)
        \sum_{k=0}^{\infty} b_{k}^{(n)}z^{-k},\ b_{0}^{(n)}=1
        \end{equation}
  \end{proposition}

\begin{remark}
For $\Pt$, Equation
(\ref{boutroux}) leads to
\[
f_0{g_1}^{2} {{\rm d}^{2}\over{\rm d} z^{2}}u
+\bigl(2 f_1g_1+f_0 g_2\bigr){{\rm d}\over{\rm d} z}u+f_2 u=
2 f_0^{3}u^{3}+x f_0 u+\alpha_{1},\]
where we are using the notation above, i.e.
$f_m=\frac{d^{m}f}{dx^{m}}$ and $g_m=\frac{d^{m}g}{dx^{m}}$.
A maximal dominant balance is given by
\begin{align*}
    &f_0{g_1}^{2}= f_0^{3}=xf_0\\
    \Rightarrow& f_0=x^{1/2},\ g_1=x^{1/2}.
    \end{align*}
The result:
\begin{equation}
    \label{p2trans}
    V(x)=x^{1/2}u(z),\quad z=2 x^{3/2}/3
    \end{equation}
leads to
\begin{equation}
    \label{p2boutroux}
    {{\rm d}^{2}\over{\rm d} z^{2}}u=2u^{3}+u+\frac 1z\left(\frac 23
    \alpha_{1}-{{\rm d}\over{\rm d} z}u\right)+\frac{u}{9z^{2}}
    \end{equation}
The new variables $u$, $z$ in (\ref{p2trans}) were first given by
Boutroux \cite{boutroux:I}.
\end{remark}

\begin{proof}
Let us first prove the following formula:
\begin{equation}
{{\rm d}^p\over{\rm d} x^p}V(x)=
\sum_{l=0}^p\tilde a_{p,l}z^{l+1-{2n\over 2n+1}(p+1)}
{{\rm d}^l\over{\rm d} z^l} u(z),
\label{bou0}\end{equation}
for some constants $\tilde a_{p,l}$,
$\tilde a_{p,p}=\left({2n+1\over 2n}\right)^{p+1\over 2n+1}$.

First let us show by induction that
\begin{equation}
V_p=\sum_{l=0}^p{\mathcal U}_{p,l}(g_1,g_2,\dots,g_p)
{{\rm d}^l\over{\rm d} z^l} u(z)
\label{bou1}\end{equation}
where for $l=1,\dots,p$,
$$
{\mathcal U}_{p,l}(g_1,g_2,\dots,g_p)=
\underset{{\overset{
\scriptstyle n_1+2n_2+\dots +p n_p=p+1}
{\scriptstyle n_1+n_2+\dots +n_p=l+1}}}
{{\boldsymbol{\Sigma}}}
c^{k,l}_{n_1 n_2 \dots n_p}
{g_1}^{n_1}\dots{g_p}^{n_p},
$$
for some constant coefficients $c^{k,l}_{n_1 n_2 \dots n_p}$, and
$$
{\mathcal U}_{p,0}= f_{p}= g_{p+1},\qquad
{\mathcal U}_{p,p}= f_0(g_1)^{p}= (g_1)^{p+1}.
$$

The formula (\ref{bou1}) is obvious for $p=0$, with
${\mathcal U}_{0,0}=f_0=g_1$. Suppose (\ref{bou1}) is valid for some $p\geq 1$
let us prove it for $p+1$.
$$
\begin{array}{ll}
V_{p+1}=&{{\rm d}\over{\rm d}x}V_p=
{{\rm d}\over{\rm d}x}
\sum_{l=0}^p{\mathcal U}_{p,l}(g_1,g_2,\dots,g_p)
{{\rm d}^l\over{\rm d} z^l} u(z)=\\
&= \sum_{l=0}^p \left({{\rm d}\over{\rm d}x}\left({\mathcal U}_{p,l}  \right)
{{\rm d}^l\over{\rm d} z^l} u(z) + {\mathcal U}_{p,l} g_1
{{\rm d}^{l+1}\over{\rm d} z^{l+1}} u(z)\right)\\
\end{array}
$$
so that we obtain (\ref{bou1}) for $p+1$ with
$$
{\mathcal U}_{p+1,l}={{\rm d}\over{\rm d}x}{\mathcal U}_{p,l}+
g_1 {\mathcal U}_{p,l-1},\qquad \hbox{for}\, l=1,\dots,p
$$
and
$$
{\mathcal U}_{p+1,0}={{\rm d}\over{\rm d}x}{\mathcal U}_{p,0},\qquad
{\mathcal U}_{p+1,p+1}= g_1{\mathcal U}_{p,p}
$$
Let us express the polynomials ${\mathcal U}_{p,l}$ as functions of $z$:
$$
\begin{array}{ll}
{g_1}^{n_1}\dots{g_p}^{n_p}\propto&
x^{\sum_{s=1}^p\left({2n+1\over 2n}-s\right)n_s} \propto
x^{{2n+1\over 2n}(n_1+\dots+n_p)-(n1+2 n_2+\dots +p n_p)}\propto\\
&\propto z^{(l+1)-{2n\over 2n+1}(p+1)}\end{array}
$$
so that we obtain formula (\ref{bou0}), for $l=1,\dots,p-1$ and
$$
{\mathcal U}_{p,0}\propto z^{1-{2n\over 2n+1}(p+1)},\qquad
{\mathcal U}_{p,p}=\left({2n+1\over 2n}z\right)^{p+1\over 2n+1}.
$$
This shows formula (\ref{bou0}).

Of course formula (\ref{bou0}) gives
$$
V_{2n}= \sum_{l=0}^{2n}\tilde a_{2n,l}z^{l+1-2n}
{{\rm d}^l u\over{\rm d} z^l}=
{2n+1\over 2n} z {{\rm d}^{2n}u\over{\rm d} z^{2n}} +
\sum_{l=0}^{2n-1}\tilde a_{2n,l}z^{l+1-2n}
{{\rm d}^l u \over{\rm d} z^l}.
$$
Let us now compute the polynomial $P_{2n-1}$ as a function of $z$, $u$
and its derivatives. Recall formlula (\ref{polyn}):
$$
P_{2n-1}=
\underset{{\overset{
\scriptstyle \langle{\bf k}\rangle=2n+1}{\scriptstyle k_0\leq 2n-1}}}
{{\boldsymbol{\Sigma}}}
b_{k_0,\dots,k_{2n-2}} V_0^{k_0} V_1^{k_1}\dots
V_{2n-2}^{k_{2n-2}}.
$$
By (\ref{bou0}), we have
$$
\begin{array}{ll}
V_p^{k_{p}} =& \sum_{s_0^p+\dots+s_p^p=k_p}
\prod^p_{l=0}\hat a_{s_0^p,\dots,s_p^p}
\left(z^{(l+1)-{2n\over 2n+1}(p+1)}
{{\rm d}^l\over{\rm d} z^l} u(z)\right)^{s_l^p}=\\
&=z^{k_p-{2n\over 2n+1}(p+1)k_p} \sum_{s_0^p+\dots+s_p^p=k_p}
\hat a_{s_0^p,\dots,s_p^p} \prod^p_{l=0}
\left(z^{l}{{\rm d}^l\over{\rm d} z^l} u(z)\right)^{s_l^p}=\\
&= z^{k_p-{2n\over 2n+1}(p+1)k_p} \sum_{s_0^p+\dots+s_p^p=k_p}
z^{\sum_{l=0}^p l s_l^p}
\hat a_{s_0^p,\dots,s_p^p} \prod^p_{l=0}
\left({{\rm d}^l\over{\rm d} z^l} u(z)\right)^{s_l^p}.\\
\end{array}
$$
However, for $p\not=0$,
$\sum_{l=0}^p l s_l^p< p \sum_{l=0}^p s_l^p=p k_p$, for $p\neq1$ and
$\sum_{l=0}^1 l s_l^1\leq k_1$, we have that
$$
V_p^{k_{p}} = z^{k_p-{2n\over 2n+1}(p+1)k_p}
\sum_{s_0^p+\dots+s_p^p=k_p}
h_{s_0^p,\dots,s_p^p}(z) \prod^p_{l=0}
\left({{\rm d}^l\over{\rm d} z^l} u(z)\right)^{s_l^p}
$$
where $h_{s_0^p,\dots,s_p^p}(z)$ are polynomials in $z$ of degree
$< p k_p$ for $p\neq1$ and of degree $\leq k_1$ for $p=1$. Now to compute
$P_{2n-1}$ as a function of $z$ let us first observe that
$$
\prod_{p=0}^{2n-2}
z^{k_p-{2n\over 2n+1}(p+1)k_p}= z^{k_0+\dots+k_{2n-2}-2n}.
$$
Thus $V_0^{k_0} V_1^{k_1}\dots V_{2n-2}^{k_{2n-2}}
z^{2n-(k_0+\dots+k_{2n-2})}$ has the form
$$
\prod_{p=0}^{2n-2}\sum_{s_0^p+\dots+s_p^p=k_p}
h_{s_0^p,\dots,s_p^p}(z) \prod^p_{l=0}
\left({{\rm d}^l\over{\rm d} z^l} u(z)\right)^{s_l^p},
$$
where the $l$-th derivative of $u$, ${{\rm d}^l\over{\rm d} z^l}u(z)$,
appears for every $p$ as
$\left({{\rm d}^l\over{\rm d} z^l} u(z)\right)^{s_l^p}$. Assuming that
$s_l^p=0$ for all $l>p$ we can write
$$
\prod_{p=0}^{2n-2}\sum_{s_0^p+\dots+s_{2n-2}^p=k_p}
h_{s_0^p,\dots,s_{2n-2}^p}(z) \prod^{2n-2}_{l=0}
\left({{\rm d}^l\over{\rm d} z^l} u(z)\right)^{s_l^p},
$$
Let us
collect together all the derivatives of $u$ of the same order, we
obtain that the $l$-th derivative of $u$ appears with the power
$m_l=\sum_{p=0}^{2n-2}s_l^p$. Thus $\sum_{l=0}^{2n-2}m_l=
\sum_{l=0}^{2n-2}\sum_{p=0}^{2n-2}s_l^p=\sum_{p=0}^{2n-2} k_p$. Now
the sum $\sum_{p=0}^{2n-2} k_p$ is maximal for  $k_0=2n-1$, $k_1=1$
and $k_p=0$ for all $p\neq 0,1$ because $\sum_{p=0}^{2n-2}(p+1) k_p=
2n+1$ and $k_0\leq 2n-1$. As a consequence we have
$0\leq m_l\leq \sum_{p=0}^{2n-2} k_p\leq 2n$.
We then obtain
$$
V_0^{k_0} V_1^{k_1}\dots V_{2n-2}^{k_{2n-2}}
z^{2n-\sum_0^{2n-2} k_p}=
\sum_{m_0,\dots,m_{2n-2}=0}^{2n}
\tilde h_{m_0,\dots,m_{2n-2}}(z)\prod_{l=0}^{2n-2}
\left({{\rm d}^l\over{\rm d} z^l} u(z)\right)^{m_l}
$$
for some polynomials $\tilde h_{m_0,\dots,m_{2n-2}}(z)$ of degree
$<\sum_{p=0}^{2n-2}p k_p$. In fact $V_1^{k_1}$ contributes with a
polynomial of degree $k_1$ and for $p\neq1$ each $V_p^{k_p}$ contributes
with a polynomial of degree $<p k_p$.
Since $\sum_{p=0}^{2n-2}(p+1)k_p=2n+1$,
for $k_1\neq0$ ($2k_1$ is even)
at least some other $k_p$ must be non-zero giving total
degree $<\sum_{p=0}^{2n-2}p k_p$.
Since $a_{m_0,\dots,m_{2n-2}}(z):=
z^{k_0+\dots+k_{2n-2}}\tilde h_{m_0,\dots,m_{2n-2}}(z)$ are
polynomials of degree $<\sum_{p=0}^{2n-2}(p+1) k_p=2n+1$, we
obtain
$$
P_{2n-1}=z^{-2n} \sum_{m_0,\dots,m_{2n-2}=0}^{2n}
a_{m_0,\dots,m_{2n-2}}(z)\prod_{l=0}^{2n-2}
\left({{\rm d}^l\over{\rm d} z^l} u(z)\right)^{m_l},
$$
where $a_{m_0,\dots,m_{2n-2}}$ have the desired properties.
This concludes the proof of (\ref{bou4}).

From (\ref{bou4}) it is straightforward to compute the leading terms
of (\ref{series1}) and (\ref{series2}) for the solutions $u(z)$
in the case when $d^lu/dz^l\ll u$, for $|z|\gg1$.
Moreover, rewriting the equation as
\begin{equation*}
\begin{array}{ll}
{2n\over 2n+1} u(z)\, +&
\beta_n u^{2n+1}+{2n\over 2n+1}{\alpha\over z}=
-{{\rm d}^{2n}\over{\rm d} z^{2n}} u
-\sum_{l=0}^{2n-1} a_{2n,l}z^{l-2n}
{{\rm d}^l\over{\rm d} z^l} u(z)+\\
&+z^{-2n-1}\sum_{m_0,\dots,m_{2n-2}=0}^{2n-1}
a_{m_0,\dots,m_{2n-2}}(z)\prod_{l=0}^{2n-2}
\left({{\rm d}^l\over{\rm d} z^l} u(z)\right)^{m_l},\\
\end{array}
\end{equation*}
and iterating on partial sums of (\ref{series1}) or (\ref{series2}),
we are led to formal solutions of the desired form.
\end{proof}

We now use Wasow's Theorem 12.1 of \cite{Was} to show the existence of
true solutions on proper sub--sectors of angular width $<\pi$. We can
write our equation in vector form
$$
\ddz Y = F(Y,z)
$$
where $Y$ is a column vector of $2n$ components
$y_j={{\rm d}^{j-1}\over{\rm d} z^{j-1}} u(z)$ and $F(Y,z)$ is the
vector function of $2n$ components $f_j$ given by $f_j=y_{j+1}$ for
$j=1,\dots,2n-1$ and
$$
\begin{array}{ll}
f_{2n}=&-\sum_{l=0}^{2n-1} a_{2n,l}z^{l-2n} y_{l+1}+ y_1 +
\beta_n y_1^{2n+1}+{\alpha\over z}+\\
&+z^{-2n-1}
\overset{2n}{\underset{m_0,\dots,m_{2n-2}=0}
{{\boldsymbol{\Sigma}}}}
a_{m_0,\dots,m_{2n-2}}(z)\prod_{l=0}^{2n-2}
y_{l+1}^{m_l}\\
\end{array}
$$
Such a function is analytic in $z$ for $|z|>z_0$, for some value
$z_0$. Let us compute the Jacobian at the point $Y=0$. Obviously, for
all $j=1,\dots,2n-1$ we have ${\partial f_j\over\partial y_i}=0$ for
all $i\neq j+1$ and ${\partial f_j\over\partial y_{j+1}}=1$.
Now let us deal with
${\partial f_{2n}\over\partial y_i}$. Observe that the term
$z^{-2n-1}\overset{2n}{\underset{m_0,\dots,m_{2n-2}=0}
{{\boldsymbol{\Sigma}}}}
a_{m_0,\dots,m_{2n-2}}(z)\prod_{l=0}^{2n-2}y_{l+1}^{m_l}$
has at least quadratic term in $y_1,\dots,y_{2n}$ because
$\sum_{l=0}^{2n-2}m_l=\sum_{p=0}^{2n-2} k_p\geq 2$, thus its
contribution to the  Jacobian at the point $Y=0$ is $0$. We obtain
then ${\partial f_{2n}\over\partial y_1}=1$ and
${\partial f_{2n}\over\partial y_i}=0$ for $i>1$.

This shows that all hypotheses of Theorem 12.1 in \cite{Was} are fulfilled and
thus the existence of true solutions on proper sub--sectors
of angular width $<\pi$ is proved.

This concludes the prood of the first and
second statements of our Theorem \ref{main}. To prove the third statement
we need the following

\begin{theorem} Consider the linear system
\begin{equation}
\eps^h \ddz Y = A(z,\eps) Y,
\label{w1}\end{equation}
where $Y$ is a column vector of $m$ components, $z\in\Complex$, $h$
is a positive integer and $A(z,\eps)$ is a $m\times m$ matrix function
admitting asymptotic expansion of the form
\begin{equation}
A(z,\eps)\sim \sum_0^\infty A_r(z)\eps^r,\qquad \eps\to 0,
\label{w2}\end{equation}
that is uniformly valid for
$\eps\in\Sigma\cap\{\eps|0<\eps\leq\eps_0\}$,
$\Sigma$ being some sector of the $\eps$-plane, with coefficients
$A_r(z)$ holomorphic in a ball, $z\in B(z_0,a)$. Denote the eigenvalues
of $A_0(z)$ by $\mu_1(z),\dots,\mu_m(z)$.
If $\mu_1(z_0),\dots,\mu_n(z_0)$ are pairwise distinct, then the system
(\ref{w1}) admits a fundamental solution $Y(z)$ of the form
\begin{equation}
Y(z)=\hat Y(z,\eps) \exp(Q(z,\eps)),
\label{w3}\end{equation}
where $\hat Y(z,\eps)$ admits asymptotic expansion of the form
\begin{equation}
\hat Y(z,\eps)\sim \sum_0^\infty\hat Y_r(z)\eps^r,\qquad \eps\to 0,
\label{w4}\end{equation}
that is uniformly valid for
$\eps\in\hat\Sigma\cap\{\eps|0<\eps\leq\hat\eps_0\}$, $\hat\eps_0<\eps$,
$\hat\Sigma$ being a
small enough sub-sector of $\Sigma$ centered around any arbitrary ray
${\rm arg}(\eps)=\alpha$, for any $\alpha$ that is not an odd multiple of
${\pi\over2}$, with coefficients $\hat Y_r(z)$ holomorphic
in $z\in B(z_0,\hat a)\subset B(z_0,a)$, $\hat a>a$, and with
$\det\hat Y_0(z)\neq 0$. The matrix $Q(z,\eps)$ is diagonal of the form
\begin{equation}
Q(z,\eps)=\sum_1^h\hat Q_r(z)\eps^{-r},\
\label{w5}\end{equation}
where $Q_h(z)=\rm{diag}\left(\int_{z_0}^z\mu_1(t)\rm{d}t,\dots,
\int_{z_0}^z\mu_m(t)\rm{d}t\right)$.
\label{Was}\end{theorem}

\proof It can be found in \cite{Was,Sib}.

Let us consider two adjacent sectors $\tilde{\mathcal S}_n$ and
$\hat{\mathcal S}_n$ and two solutions, ${\tilde V}(x)$ and
${\hat V}(x)$ of (\ref{n1}) having the same asymptotic behaviour
(\ref{n2.5}) or (\ref{n3.5}) in the given sectors $\tilde{\mathcal S}_n$
and $\hat{\mathcal S}_n$ respectively.
For every fixed $j=0,\dots,2n$, we take, in the case of asymptotic
behavior (\ref{n2.5}), the sectors
$$
\begin{array}{ccc}
     \tilde{\mathcal S}_n&=&\left\{x\in\Complex\bigm| |x|>|x_{0}|,
   {(2j-1)\pi\over 2n+1}+\epsilon
<|\arg(x)|<{(2j+1)\pi\over 2n+1}+\epsilon\right\},\\
  \hat{\mathcal S}_n&=&\left\{x\in\Complex\bigm| |x|>|x_{0}|,
  {(2 j+1) \pi\over 2n+1}-\epsilon<
|\arg(x)|<{(2j+3)\pi\over 2n+1}-\epsilon \right\},\\
  \end{array}
$$
and in the case of asymptotic behavior (\ref{n3.5}) we take the sectors
$$
\begin{array}{c}
     \tilde{\mathcal S}_n=\left\{x\in\Complex\bigm| |x|>|x_{0}|,
   \pi+{(2j-1)\pi\over 2n+1}+\epsilon
<|\arg(x)|<\pi+{(2j+1)\pi\over 2n+1}+\epsilon\right\},\\
  \hat{\mathcal S}_n=\left\{x\in\Complex\bigm| |x|>|x_{0}|,
  \pi+{(2j+1)\pi\over 2n+1}-\epsilon<
|\arg(x)|<\pi+{(2j+3)\pi\over 2n+1}-\epsilon \right\},\\
  \end{array}
$$
for $\epsilon>0$ small enough. In both cases
these sectors intersect in a small sector ${\mathcal S}_\epsilon$
of opening $<2\epsilon$ containing in the first case the line
$$
{\rm arg}(x) =  {2j+1\over 2n+1}\pi,
$$
and, in the second case, the line
$$
{\rm arg}(x) = \pi + {2j+1\over 2n+1}\pi.
$$
We want to show that in this small sector ${\mathcal S}_\epsilon$ the
two solutions ${\tilde V}(x)$ and ${\hat V}(x)$ coincide. If this is the case,
we can then analytically extend ${\tilde V}(x)$ to the sector
$\Sigma_n=\tilde{\mathcal S}_n\cup \hat{\mathcal S}_n$, this extension is
unique and the third statement of our Theorem \ref{main} is proved.

Suppose by contradiction that the  two solutions, ${\tilde V}(x)$
and ${\hat V}(x)$ of (\ref{n1}) having the same asymptotic behaviour
(\ref{n2.5}) or (\ref{n3.5}) differ in the given sector
${\mathcal S}_\epsilon$. Then
their difference $W(x):={\tilde V}(x)-{\hat V}(x)$ and all its derivatives
vanish asymptotically as $x\to\infty$ in ${\mathcal S}_\epsilon$. Suppose
that $W(x)\neq 0$ for some value $x$.

Let us consider the differential equation satisfied by $W(x)$
\begin{equation}
W_{2 n}= P_{2n-1}({\tilde V}_0,\dots,{\tilde V}_{2n-2})-
P_{2n-1}({\hat V}_0,\dots,{\hat V}_{2n-2}) + x W_0+
\beta_n ({\tilde V}_0^{2n+1}-{\hat V}_0^{2n+1}),
\label{diff1}\end{equation}
where ${\tilde V}_m=\frac{d^{m}\tilde V}{dx^{m}}$,
${\hat V}_m=\frac{d^{m}\hat V}{dx^{m}}$ and
$W_m=\frac{d^{m}W}{dx^{m}}$ as above.
Now
$$
{\tilde V}_0^{2n+1}-{\hat V}_0^{2n+1}
=\sum_{k=0}^{2n}{\tilde V}_0^{2n-k}{\hat V}_0^{k}\,W_0
$$
and
$$
\begin{array}{l}
P_{2n-1}({\tilde V}_0,\dots,{\tilde V}_{2n-2})-
P_{2n-1}({\hat V}_0,\dots,{\hat V}_{2n-2})= \\
\underset{{\overset{
\scriptstyle \langle{\bf k}\rangle=2n+1}{\scriptstyle k_0\leq 2n-2}}}
{{\boldsymbol{\Sigma}}}
b_{k_0,\dots,k_{2n-2}}\left({\tilde V}_0^{k_0} {\tilde V_1}^{k_1}\dots
{{\tilde V}_{2n-2}}^{k_{2n-2}}-{\hat V}_0^{k_0} {{\hat V_1}}^{k_1}\dots
{{\hat V}_{2n-2}}^{k_{2n-2}}\right)\\
=Q_{2n-1} W_0 +  Q_{2n-2} W_1 + \dots
Q_{1} W_{2n-2}\\
\end{array}
$$
where
$$
Q_{2n-p-1}=\tilde b^{(p)}_{k_0,\dots,k_{2n-2}}
{\hat V}_0^{k_0} {\hat V}_1^{k_1}\dots
{\hat V}_{p-1}^{k_{p-1}} \left(\sum_{l=0}^{k_p-1}
{\tilde V}_{p}^{k_{p}-1-l}{\hat V}_{p}^{l}\right)
{\tilde V}_{p+1}^{k_{p+1}}\dots {\tilde V}_{2n-2}^{k_{2n-2}},
$$
for some constants $\tilde b^{(p)}_{k_0,\dots,k_{2n-2}}$.
So we have
\begin{equation}
W_{2 n}= Q_{2n-1} W_0 +  Q_{2n-2} W_1 + \dots + Q_{1} W_{2n-2} + x W_0 +
\beta_n \sum_{k=0}^{2n}{\tilde V}_0^{2n-k} {\hat V}_0^{k}\,W_0.
\end{equation}
We want to study the asymptotic behaviour of the non-zero solutions
$W(x)$ of this linear differential equation as $x\to\infty$ in the
sector $\Sigma_n$.

We first deal with the case (\ref{n2.5}). Let us perform a variable
rescaling $\eps^{2n} x=z$, such that $z\in B(z_0,a)$ for some $a>0$
and $\eps\to 0$, in the sector $\Sigma_n$. From (\ref{n2.5}) we obtain
$$
\tilde V\left({z\over\eps^{2n}}\right),\,
\hat V\left({z\over\eps^{2n}}\right) \sim
\eps^{-1} a_0 z^{1\over 2n} + {\mathcal O}(\eps^{2n})
$$
so that the polynomials $Q_{2n-p-1}$ are all rescaled as
$Q_{2n-p-1}\to \eps^{2n+2} q_{2n-p-1}$, for some polynomials
$q_{2n-p-1}(z)$. In fact
$V_p^{k_p} \sim \epsilon^{k_p(2np-1)} {{\rm d}^p v\over{\rm d}z^p}$
and since
$$
\begin{array}{l}
k_0 (-1)+k_1(2n-1)+\ldots+k_{p-1}(2n(p-1)-1)+(k_p-1)(2np-1)+\\
+k_{p+1}(2n(p+1)-1)+ \ldots+ k_{2n-2}(2n(2n-2)-1)=\\
=\sum_{l=0}^{2n-2} k_l(2nl-1)-(2n p-1)=
2n\left(\sum_{l=0}^{2n-2} k_l l -p\right)-\sum_{l=0}^{2n-2} k_l +1\\
\geq 2n\left(\sum_{l=0}^{2n-2} k_l l -p\right)-2n+2\geq 2n-2,\\
\end{array}
$$
because $\sum_{l=0}^{2n-2} k_l l\geq \sum_{l=0}^{2n-2} k_l+1$, we have
$$
{\hat V}_0^{k_0} {\hat V}_1^{k_1}\dots {\hat V}_{p-1}^{k_{p-1}}
\left(\sum_{l=0}^{k_p-1} {\tilde V}_{p}^{k_{p}-1-l}
{\hat V}_{p}^{l}\right)
{\tilde V}_{p+1}^{k_{p+1}}\dots {\tilde V}_{2n-2}^{k_{2n-2}}=
{\mathcal O}\left(\eps^{2(n+1)}\right).
$$
Analogously the polynomial $\sum_{k=0}^{2n}{\tilde V}^{2n-k}{\hat V}^{k}$
can be expanded as $(2n+1){a_0^{2n}\over\eps^{2n}} z+
{\mathcal O}(\eps)$. The differential equation then becomes:
\begin{equation}\begin{array}{ll}
\eps^{4 n^2} {{\rm d}^{2n}\over{\rm d} z^{2n}}W=& \eps^{2n+2}
\sum_{k=0}^{2n-2} q_{2n-k-1}(z){{\rm d}^{k}\over{\rm d} z^{k}}W+\\
&+ {z\over\eps^{2n}} W - {2n\over\eps^{2n}} z W +{\mathcal O}(\eps) W.\\
\end{array}\end{equation}
Such equation can be put into system form by
$$
Y_1=  W(z),\qquad Y_{k+1}= \eps^{2n+1} {{\rm d}\over{\rm d} z}Y_k,
\qquad k=1,\dots,2n-1.
$$
In this way we obtain a system of the form (\ref{w1}) with $h=2n+1$ and
leading matrix
\begin{equation}
A_0(z)=\left(\begin{array}{ccccc}
0&1&0&\dots&0\\
0&0&1&\dots&0\\
\dots&\dots&\dots&\dots&\dots\\
0&\dots&\dots&0&1\\
f(z)&0&\dots&\dots&0\\
\end{array}\right),
\label{mata0}\end{equation}
with $f(z):=(1-2n)z$.
The eigenvalues of $A_0$ are
$\mu_k(z)=(2n-1)^{1\over 2n}z^{1\over 2n}\exp({k\pi i\over n})$,
$k=1,\dots,2n$. Applying Theorem \ref{Was}, we conclude that there exists
a fundamental solution $Y(z)$ of our system of the form (\ref{w1}) with
asymptotic behaviour
$$
Y(z)\sim \sum_0^\infty\hat Y_r(z)\eps^r
\exp\left(
{1\over\eps^h}{\rm diag}\left(\nu_1(z),\dots,\nu_{2n}(z)\right)+
{\mathcal O}(\eps^{1-h})
\right)
$$
where $\nu_k(z)={2n(2n-1)^{1\over 2n}\over 1+2n}
\exp({i\pi k\over n})z^{1+2n\over 2n}$ and ${\rm det}(\hat Y_0(z))\neq 0$.
Such a fundamental solution $Y(z)$ exists for every sub-sector of
${\mathcal S}_\epsilon$ centered around the line (\ref{n2.66}).
In fact such lines are never odd
multiples of ${\pi\over2}$. This shows that in any sub-sector of
${\mathcal S}_\epsilon$ centered around the line (\ref{n2.66})
each solution $W(x)$ of (\ref{diff1}) has asymptotic
behaviour with leading term
$$
{W}(x)\sim {W}_0\,{\rm diag}\left(
\exp\left({2n(2n-1)^{1\over 2n}\over 1+2n} \exp({i\pi k\over n})
x^{1+2n\over 2n}\right), \, k=1,\dots, 2n\right)
$$
where ${W}_0$ is a nonsingular constant matrix. This leads to a
contradiction because the asymptotic behaviour of
each solution $W(x)$ of (\ref{diff1}) or of at least one of its
derivatives is oscillatory along those lines, and does not vanish
asymptotically as assumed at the beginning.

We now deal with the case (\ref{n3.5}). Let us perform a variable
rescaling $\tau x=z$, such that $z\in B(z_0,a)$ and $\tau\to 0$, in the
sector $\Sigma_n$. Then $\tilde V$ and $\hat V$ are accordingly
rescaled to $\tau\tilde v$ and $\tau\hat v$ respectively.
As a consequence
$$
\begin{array}{ll}
Q_{2n-p-1}W_p\to &
\tau^{k_0+\dots+p k_{p-1}+(p+1)(k_p-1)+(p+2) k_{p+1}+\dots+(2n-1)k_{2n-2}}
\tau^{p+1}{{\rm d}^p W\over{\rm d} z^p}=\\
&=\tau^{\langle{\bf k}\rangle-(p+1)}
\tau^{p+1}{{\rm d}^p W\over{\rm d}} z^p=
\tau^{2n+1} q_{2n-p-1}{{\rm d}^p W\over{\rm d} z^p} \\
\end{array}
$$
for some polynomial $q_{2n-p-1}$. The differential
equation then becomes:
\begin{equation}
\begin{array}{ll}
\tau^{2 n} {{\rm d}^{2n}\over{\rm d} z^{2n}}W= &
Q_{2n} W + \tau Q_{2n-1}{{\rm d}\over{\rm d} z} W + \dots +
\tau^{2n-1}Q_{1} {{\rm d}^{2n-1}\over{\rm d} z^{2n-1}}W
+\\ & + {z\over\tau} W +
\beta_n \tau^{2n}({\tilde v}^{2n}+\dots+{\hat v}^{2n})W.\\
\end{array}\end{equation}
By the rescaling
$$
Y_1(z):=W(z),\quad Y_{k+1}:=\eps^h {{\rm d}\over{\rm d} z}Y_k(z),
\qquad k=1,\dots,2n-1,
$$
where $\eps=\tau^{1\over 2 n}$ and $h=2n+1$ we obtain again a system of the
form (\ref{w1}) with leading matrix $A_0(z)$ of the same form (\ref{mata0})
with $f(z)=z$. The eigenvalues of $A_0(z)$ are
$\mu_k(z)=z^{1\over 2n}\exp\left({i k\pi\over n}\right)
\exp\left({i\pi\over2n}\right)$, $k=1,\dots,2 n$. Applying Theorem
\ref{Was}, we conclude that there exists a fundamental solution
$Y(z)$ of our system
of the form (\ref{w1}) with asymptotic behaviour
$$
Y(z)\sim \sum_0^\infty\hat Y_r(z)\eps^r
\exp\left(
{1\over\eps^h}{\rm diag}\left(\nu_1(z),\dots,\nu_{2n}(z)\right)
\right)
$$
where $\nu_k(z)={2n\over 1+2n}\exp\left({i\pi\over2n}\right)
\exp({i\pi k\over n})z^{1+2n\over 2n}$ and ${\rm det}(\hat Y_0(z))\neq 0$.
Such solution $Y(z)$ exists for every sub-sector of
${\mathcal S}_\epsilon$ centered around the line (\ref{n3.66}).
In fact such lines are never odd multiples of
${\pi\over2}$.  This shows that in any sub-sector of
${\mathcal S}_\epsilon$ centered around the line (\ref{n3.66})
each solution $W(x)$ of (\ref{diff1}) has asymptotic
behaviour with leading term
$$
W(x)\sim {W}_0\,{\rm diag}\left(
\exp\left({2n\over 1+2n}\exp\left({i\pi\over2n}\right)
\exp\left({i\pi k\over n}\right)x^{1+2n\over 2n}\right)
\, k=1,\dots, 2n\right)
$$
where ${W}_0$ is a non-singular constant matrix.
This leads to a contradiction because the asymptotic behaviour of
each solution $W(x)$ of (\ref{diff1}) or of at least one of its
derivatives is oscillatory along those lines, and not vanishing
asymptotically
as assumed at the beginning.


\begin{thebibliography}{1}

\bibitem{AS}
M.~J. Ablowitz and H. Segur.
\newblock Exact linearization of a Painlev\'e transcendent.
\newblock {\em Phys. Rev. Lett.}, 38:1103--1106, 1977.

\bibitem{Ai}
H. Airault.
\newblock Rational solutions of Painlev\'e equations.
\newblock {\em Studies in Applied Mathematics}, 61:31--53, 1979.

\bibitem{boutroux:I}
P.~Boutroux.
\newblock Recherches sur les transcendantes de {M. Painlev\'e} et l'\'etude
  asymptotique des \'equations diff\'erentielles du second ordre.
\newblock {\em Ann. \'Ecole Norm.}, 30:265--375, 1913.

\bibitem{boutroux:II}
P.~Boutroux.
\newblock Recherches sur les transcendantes de {M. Painlev\'e} et l'\'etude
  asymptotique des \'equations diff\'erentielles du second ordre.
\newblock {\em Ann. \'Ecole Norm.}, 31:99--159, 1914.

\bibitem{CJM}
P.~A. Clarkson and N.~Joshi and M.~Mazzocco.
\newblock The Lax pair for the \text{mKdV} hierarchy.
\newblock {\em preprint}, 2002.


\bibitem{cjp}
P.~A. Clarkson and N.~Joshi and A. Pickering.
\newblock B\"acklund transformations for the second {Painlev\'e}
hierarchy: a modified truncation approach.
\newblock {\em Inverse Problems}, 15:175--187, 1999.




\bibitem{FN}
H. Flaschka and A.~C. Newell.
\newblock Monodromy and Spectrum Preserving Deformations I.
\newblock {\em Comm. Math. Phys.}, 76:65--116, 1980.

\bibitem{Hab1}
R. Haberman.
\newblock The modulated phase shift for weakly dissipated nonlinear
oscillatory waves of the Korteweg-de Vries type.
\newblock {\em Stud. Appl. Math.}, 78,  no. 1, 73--90 1988.


\bibitem{Hab2}
R. Haberman.
\newblock Nonlinear transition layers---the second Painlev\'e
transcendent.
\newblock {\em Stud. Appl. Math.}, 57, no. 3, 247--270, 1977.


\bibitem{njmdk:conn1}
N.~Joshi and M.~D. Kruskal.
\newblock The {Painlev\'e} connection problem: an asymptotic approach {I}.
\newblock {\em Stud. Appl. Math.}, 86:315--376, 1992.

\bibitem{lax}
P.~D. Lax.
\newblock Almost periodic solutions of the KdV equation.
\newblock{\em SIAM review}, 18, no.3:351--375, 1976.

\bibitem{Li}
R.~C. Littlewood.
\newblock Hyperelliptic asymptotics of \Pa--type equations.
\newblock{\em Nonlinearity}, 12, 1629--1641, 1999.

\bibitem{Pain}
P. Painlev\'e.
\newblock Sur les Equations Differentielles du Second Ordre et
d'Ordre Superieur, dont l'Interable Generale est Uniforme.
\newblock {\em Acta Math.}, 25:1--86, 1902.

\bibitem{Sib}
Y. Sibuya.
\newblock Linear Differential Equations in the Complex Domain: Problems
of Analytic Continuation.
\newblock {\em  AMS TMM}, 82, 1990.

\bibitem{Um}
H. Umemura.
\newblock Second Proof of the Irreducibility of the First
Differential Equation of Painlev\'e.
\newblock {\em Nagoya Math. J.}, 117:125--171, 1990.

\bibitem{UW}
H. Umemura and H. Watanabe.
\newblock Solutions of the Second and Fourth Painlev\'e Equations. I.
\newblock {\em Nagoya Math. J.}, 148:151--198, 1997.

\bibitem{Was}
W. Wasow.
\newblock Asymptotic Expansions for Ordinary differential Equations.
\newblock {\em Pure and Applied Mathematics, John Wiley \and Sons,
Inc.}, XIV, 1965.


\end{thebibliography}
\end{document}